\documentclass{amsart}
\usepackage{amsmath, amsthm, amssymb, amscd, mathrsfs, eucal, epsfig}
\begin{document}

\newtheorem{theorem}{Theorem}[section]
\newtheorem{lemma}[theorem]{Lemma}
\newtheorem{corollary}[theorem]{Corollary}
\newtheorem{conjecture}[theorem]{Conjecture}
\newtheorem{question}[theorem]{Question}
\newtheorem{problem}[theorem]{Problem}
\newtheorem{prop}[theorem]{Proposition}
\newtheorem*{claim}{Claim}
\newtheorem*{criterion}{Criterion}
\newtheorem*{gap_thm}{Theorem A}
\newtheorem*{improved_gap_thm}{Theorem A${}^\prime$}
\newtheorem*{accumulation_thm}{Theorem B}
\newtheorem*{mod_grp_thm}{Theorem C}
\newtheorem*{amalgam_thm}{Theorem D}
\newtheorem*{essential_cor}{Corollary E}

\theoremstyle{definition}
\newtheorem{definition}[theorem]{Definition}
\newtheorem{construction}[theorem]{Construction}
\newtheorem{notation}[theorem]{Notation}

\theoremstyle{remark}
\newtheorem{remark}[theorem]{Remark}
\newtheorem{example}[theorem]{Example}
\newtheorem{kf}[theorem]{{\bf KF}}

\def\square{\hfill${\vcenter{\vbox{\hrule height.4pt \hbox{\vrule width.4pt
height7pt \kern7pt \vrule width.4pt} \hrule height.4pt}}}$}

\def\cl{\textnormal{cl}}
\def\scl{\textnormal{scl}}
\def\area{\textnormal{area}}
\def\id{\textnormal{id}}
\def\H{\mathbb H}
\def\Z{\mathbb Z}
\def\R{\mathbb R}
\def\Q{\mathbb Q}
\def\PSL{\textnormal{PSL}}
\def\SL{\textnormal{SL}}
\def\C{\mathbb C}
\def\F{\mathcal F}
\def\Mod{\textnormal{MCG}}
\def\CC{\mathcal C}
\def\CAT{\textnormal{CAT}}
\def\til{\widetilde}
\def\length{\textnormal{length}}
\def\isom{\textnormal{Isom}}
\def\axis{\textnormal{axis}}
\def\homeo{\textnormal{Homeo}}
\def\tr{\tau}
\def\inte{\textnormal{int}}

\title{Stable commutator length in word-hyperbolic groups}
\author{Danny Calegari}
\address{Department of Mathematics \\ California Institute of Technology \\
Pasadena CA, 91125}
\email{dannyc@its.caltech.edu}
\author{Koji Fujiwara}
\address{Graduate School of Information Science \\ Tohoku University \\ Sendai, Japan}
\email{fujiwara@math.is.tohoku.ac.jp}

\begin{abstract}
In this paper we obtain uniform positive lower bounds on the stable commutator length 
of elements in word-hyperbolic groups and certain groups acting on hyperbolic spaces 
(namely the mapping class group acting on the complex of curves, and an 
amalgamated free product acting on an associated Bass-Serre tree). If $G$ is a word-hyperbolic group 
that is $\delta$-hyperbolic with respect to a symmetric generating set $S$, then there 
is a positive constant $C$ depending only on $\delta$ and on $|S|$ such that every element of 
$G$ either has a power which is conjugate to its inverse, or else the 
stable commutator length of the element is at least equal to $C$. By Bavard's theorem, these 
lower bounds on stable commutator length imply the existence of quasimorphisms 
with uniform control on the defects; however, we show how to construct 
such quasimorphisms directly.

We also prove various separation theorems on families of elements in such groups, 
constructing homogeneous 
quasimorphisms (again with uniform estimates) which are positive on 
some prescribed element while vanishing on some family of independent 
elements whose translation lengths are uniformly bounded.

Finally, we prove that the first accumulation point for 
stable commutator length in a torsion-free word-hyperbolic group is contained 
between $1/12$ and $1/2$. This gives a 
universal sense of what it means for a conjugacy class in a hyperbolic group 
to have a small stable commutator length, and can be thought of as a kind 
of ``homological Margulis lemma''.
\end{abstract}

\date{6/22/2009, Version 0.28}

\maketitle

\section{Introduction}

Let $G$ be a group, and let $[G,G]$ denote the commutator subgroup. Given $g \in [G,G]$ the 
{\em commutator length} of $g$, denoted $\cl(g)$,
is the least number of commutators in $G$ whose product is equal to $g$. The
{\em stable commutator length}, denoted $\scl(g)$, is the limit of
$\cl(g^n)/n$ as $n$ goes to infinity. If $g^n \in [G,G]$ for some least positive integer $n$,
define $\scl(g) = \scl(g^n)/n$, and define $\scl(g) = \infty$ if no power of $g$ is
contained in $[G,G]$; see \S~\ref{background_section} for precise definitions. 
For a general introduction to the theory of stable commutator length, see \cite{Calegari_scl_monograph}.

Informally, if $G = \pi_1(X)$ for some topological space $X$, commutator length is the smallest
genus of a surface in $X$ whose boundary represents a given loop in $X$, and
stable commutator length ``rationally'' measures the same quantity. At the
homological level, stable commutator length is an $L^1$ filling norm with $\Q$ coefficients, 
in the sense of Gersten 
\cite{Gersten} or Gromov \cite{Gromov_bounded}.

This paper is concerned with obtaining uniform positive lower bounds on stable commutator
length in (word-) hyperbolic groups and groups acting on ($\delta$-) hyperbolic spaces. 
Morally, ($\delta$-) hyperbolic spaces are 
those whose geometry can be efficiently probed by maps of triangles and other 
surfaces into the space, so it should not be surprising that they can be studied
effectively with stable commutator length.
However, there is a sense in which our
results {\em are} counterintuitive, which we briefly explain. 

It is a widely observed fact that in the presence
of (coarse) negative curvature, one can obtain upper bounds on the ``size'' of
a surface which is efficient in some sense. Here size is measured with respect to
some kind of norm; for instance area, Whitney $\flat$-norm, $L^1$ norm on homology
etc. Efficiency might vary from context to context (e.g. harmonic, minimal, normal)
and will depend on the way in which we measure the size of the surface. 
Our main results say that in the
presence of (coarse) negative curvature, there is a uniform (positive) {\em lower} bound on the
(homological) size of certain surfaces.

We concentrate on {\em worst-case} behavior, rather than {\em typical behavior},
and our results are frequently sharp. Other authors have studied commutator
length and its relation with negative curvature, especially Gromov \cite{Gromov_asymptotic},
\S 6.$C_2$. In our language, Gromov observes that in a word-hyperbolic group,
if $g \in [G,G]$ is not torsion and does not have a power conjugate to its inverse, then
$\scl(g)>0$ (actually, Gromov neglects to mention the second possibility).
The key innovation in our paper is that our bounds are
{\em uniform}, and depend only on macroscopic features of the group,
specifically $\delta$ and number of generators. Note that any hyperbolic group may
be made $\delta$-hyperbolic for some universally small $\delta$ (say $\delta < 10$)
just by increasing the number of generators, so our results in this sense are best
possible.

\medskip

Perhaps our most striking theorem is
Theorem B, a kind of ``spectral gap'' theorem,
which says that in a non-elementary torsion-free word-hyperbolic group,
the first accumulation point for stable commutator length as a function from
conjugacy classes to positive real numbers is between $\frac 1 {12}$ and $\frac 1 2$.
These bounds imply that there is a {\em universal} sense
of what it means for an element in a word-hyperbolic group to have a small
stable commutator length. This should perhaps be compared with Margulis' Lemma,
which says that there is a universal sense of what it means for a closed geodesic in
a hyperbolic manifold to be {\em short}. The difficulty in transporting Margulis'
Lemma to geometric group theory is that geometric notions in a group are typically
only defined up to a certain ambiguity (e.g. quasi-isometry) which obscures
small scale geometric phenomena. The advantage of working with stable commutator
length, and the power of our theorem, is precisely that it captures such small
scale phenomena. This comparison is more than superficial ---
the key to obtaining the sharp lower bound is to use estimates due to Mineyev
which reproduce, in a general $\delta$-hyperbolic space, geometric phenomena
which are strictly analogous to $C^1$ phenomena in Lie groups. 

We are able to give similar uniform lower bounds on the stable commutator length of certain elements
in two other important classes of groups: mapping class groups of surfaces (hereafter denoted $\Mod(S)$), 
and amalgamated free products. In general, neither kind of group is $\delta$-hyperbolic,
but each acts naturally on a certain $\delta$-hyperbolic space: the complex
of curves, and the Bass-Serre tree respectively. Obtaining lower bounds on
stable commutator length in mapping class groups is intimately tied to important problems in
$4$-dimensional symplectic geometry; for example, estimating the complexity of a symplectic 
$4$-manifold filling a given contact $3$-manifold, or controlling the ratios of characteristic
numbers (e.g. Euler characteristic and signature). This point of view has been
pioneered by D. Kotschick, sometimes in collaboration with H. Endo,
in a series of papers which include \cite{Kotschick_monopoles}, \cite{Endo_Kotschick:invent},
\cite{Kotschick} and \cite{Endo_Kotschick}. From another point of view, there
are relations between symplectic geometry and stable commutator length which are
more intimately connected with invariants like Hofer length and subgroup distortion;
see work of Polterovich, especially \cite{Polterovich_growth}. 

\medskip

The Bavard Duality Theorem (see Theorem~\ref{Bavard_duality_theorem}) gives a duality between
stable commutator length, and certain kinds of functions on a group, called {\em homogeneous
quasimorphisms}. A function $f:G \to \R$ is a homogeneous quasimorphism if it is
homogeneous (i.e. if it satisfies $f(g^n) = nf(g)$ for all $g \in G$ and $n \in \Z$) and if there
is a least real number $D(f)$ called the {\em defect}, for which $|f(gh)-f(g)-f(h)| \le D(f)$ for
all $g,h \in G$. The Duality Theorem says that
obtaining lower bounds for stable commutator length is equivalent to
constructing homogeneous quasimorphisms with small defects (see \S~\ref{background_section}). 
Bavard's theorem is non-constructive, and uses the Axiom of Choice (in the form of the
Hahn-Banach Theorem); however in our paper we are able to construct explicit quasimorphisms
with small defects directly. We are also able to prove various {\em separation theorems},
constructing quasimorphisms with small defects that take big values on prescribed elements
and vanish on others. This is a quantitative improvement on the kinds of separation theorems
proved or conjectured by various people in the past, and we expect it to have a
number of applications. 

\medskip

Let us stress that this paper is concerned more with developing foundations,
and understanding what (to us) seems like a fundamental algebraic/geometric inequality,
manifested in several important group theoretic contexts, than with deducing topological or other
corollaries. 

\subsection{Statement of results}

We now give a brief summary of the contents of the paper. In \S 2 we give
definitions, standardize notation, and recall some basic elements of the theory
of $\delta$-hyperbolic groups and spaces, and bounded cohomology. We spend some
time discussing Mineyev's geodesic flow space, which is a technical tool that
we find useful for establishing uniform bounds on stable commutator length directly
by means of geometry rather than using quasimorphisms (see
\S \ref{background_section} for definitions).

There are other
versions of the geodesic flow space of a word-hyperbolic group or space, and that of
Gromov \cite{Gromov} is probably best known. However, Mineyev's space has some
features that make it technically easier to use for our purposes.

\medskip

In \S 3 we give the first version of our first main result, the Gap Theorem:

\begin{gap_thm}[Gap Theorem, weak version]
Let $G$ be a word-hyperbolic group that is $\delta$-hyperbolic with respect to
a symmetric generating set $S$ with $|S|$ generators. 
Then there is a constant $C(\delta,|S|)>0$
such that for every $a \in G$, either $\scl(a) \ge C$ or else there is
some positive integer $n$ and some $b \in G$ such that $ba^{-n}b^{-1} = a^n$.
\end{gap_thm}
Note that if $ba^{-n}b^{-1} = a^n$ then $\scl(a)=0$.

It is easy to produce examples of word-hyperbolic groups which contain elements
with arbitrarily small (positive) stable commutator length, so the dependence
on $\delta$ and $|S|$ is necessary.

\medskip

In \S 4 and \S 5 we introduce so-called {\em counting quasimorphisms}, and
use them to prove the second (stronger) version
of the Gap Theorem. In the following theorem, we use the notation $\tau(a)$
to denote the {\em translation length} of an element $a \in G$, as measured in
the Cayley graph $\Gamma_S(G)$ of $G$ with respect to a fixed generating set $S$.
It is worth keeping in mind that if $a$ is not torsion, there is a positive lower bound
on $\tau(a)$ that depends only on $\delta$ and $|S|$.

\begin{improved_gap_thm}[Gap Theorem, strong version]
Let $G$ be a word-hyperbolic group that is $\delta$-hyperbolic
with respect to a symmetric generating set $S$ with $|S|$ generators.
Let $a \in G$ be a (non-torsion) element.
Assume there is no $n >0$ and no $b \in G$
with 
$ba^{-n}b^{-1}=a^n$.
Then there is a homogeneous quasimorphism $h$ on $G$
such that 
\begin{enumerate}
\item{$h(a) =1 $}
\item{The defect of $h$ is $\le C(\delta,|S|)$.}
\end{enumerate}
Moreover, let $a_i \in G$ be a collection of elements for which
$T = \sup_i \tau(a_i)$ is finite. Suppose that for all integers $n,m \ne 0$
and all elements $b \in G$ and indices $i$, there is an inequality
$$ba^nb^{-1} \ne a_i^m$$
Then there is a homogeneous quasimorphism $h$ on $G$
such that 
\begin{enumerate}
\item{$h(a) =1 $, and $h(a_i)=0$ for all $i$}
\item{The defect of $h$ is 
$\le C'(\delta, |S|, T)$}
\end{enumerate}
\end{improved_gap_thm}

The weak version of the Gap Theorem follows from the strong version. Moreover,
the first part of the strong version follows from the weak version together
with an application of Bavard's theorem, but the implication in this direction
uses the Axiom of Choice. Note that the constants
$C$ appearing in the two versions are different in either case, but related.

The second part of this strong version of the Gap Theorem can be
thought of as a kind of {\em separation theorem}, which can be
summarized in words as follows: given an element
$a$ and a finite set of elements $a_i \in G$ whose nontrivial powers are
never equal to a conjugate of a nontrivial
power of $a$, one can (explicitly) construct a homogeneous quasimorphism that vanishes on all the $a_i$
and is positive on $a$; moreover, the defect of this quasimorphism can be {\em uniformly}
controlled just in terms of the {\em translation lengths} of the $a_i$,
and in terms of $\delta$ and $|S|$.

Such separation theorems have been pursued by several authors.
For example, Polterovich and Rudnick \cite{Polterovich_Rudnick}
proved such a separation theorem for the group $\SL(2,\Z)$,
and asked if one can generalize it to hyperbolic groups; our work gives a quantitative
positive answer to their question.

\medskip

In \S 6 we state and prove the Accumulation Theorem:

\begin{accumulation_thm}[Accumulation Theorem]
Let $G$ be a torsion free non-elementary word-hyperbolic group. Then the first accumulation point
$\delta_\infty$ for stable commutator length satisfies
$$\frac 1 {12} \le \delta_\infty \le \frac 1 2$$
\end{accumulation_thm}

Here $\delta_\infty$ is the first accumulation point for the values (with multiplicity)
of the function $\scl$ on the set of conjugacy classes in $G$, thought of
as an ordered subset of $\R$. As mentioned above, we think of this theorem as a kind of
``homological Margulis Lemma''.
One can obtain similar theorems for groups with
torsion, but the statement is not so clean or universal.

\medskip

In \S 7 and \S 8 we discuss two important examples of groups acting on
$\delta$-hyperbolic spaces: the action of the mapping class group on the
complex of curves, and the action of an amalgamated free product on its
associated tree.

Our main theorem about the mapping class group is:

\begin{mod_grp_thm}[Mapping Class Theorem]
Let $S$ be a compact orientable surface of hyperbolic type and 
$\Mod(S)$ its mapping class group.
Then there is a positive integer $P$ depending on $S$ such that 
for any pseudo-Anosov element $a$, either there is
an $0 < n \le P$ and an element $b \in \Mod(S)$ with 
$ba^{-n} b^{-1}=a^n$, or else there exists a
homogeneous quasimorphism $h$ on $\Mod(S)$ such that 
$h(a) =1$ and the defect of $h$ is $\le C(S)$,
where $C(S)$ depends only on $S$.

Moreover, let $a_i \in \Mod(S)$ be a collection of elements for which
$T = \sup_i \tau(a_i)$ is finite. Suppose that for all integers $n,m \ne 0$
and all elements $b \in \Mod(S)$ and indices $i$, there is an inequality
$$ba^nb^{-1} \ne a_i^m$$
Then there is a homogeneous quasimorphism $h$ on $\Mod(S)$
such that 
\begin{enumerate}
\item{$h(a) =1 $, and $h(a_i)=0$ for all $i$}
\item{The defect of $h$ is 
$\le C'(S, T)$}
\end{enumerate}
\end{mod_grp_thm}

Note that if $b$ is not pseudo-Anosov, then $\tr(b)=0$. Moreover, for
any surface $S$, there are infinitely many conjugacy classes of elements
$b$ with $\tr(b) \le O(1/g\log(g))$ where $g = \text{genus}(S)$ (see \cite{Farb_Leininger_Margalit},
especially Theorem 1.5).

In the case of the mapping class group, one knows that not all the
quasimorphisms are of the kind promised by this theorem. Along completely
different lines, Endo and Kotschick \cite{Endo_Kotschick:invent}
managed to obtain lower bounds on stable commutator
length for elements $a \in S$ which are products of {\em positive Dehn twists in disjoint curves}. 
Such elements are reducible, and therefore fix a simplex in the complex
of curves $\CC(S)$; by contrast, our methods do not give any information about elements
with fixed points.

Note that our arguments make use of a key {\em acylindricity}
property proved by Bowditch, which is something like a bicombing in the context of $\CC(S)$.
Our lower bounds {\em do} depend on the surface $S$. There is some evidence that
this dependence is necessary: a forthcoming paper by Kotschick \cite{Kotschick_prepare}
shows that stable commutator length goes to zero under stabilizing genus by
adding handles.

Note that as in Theorem A${}^\prime$, Theorem C includes a {\em separation theorem}.
Endo and Kotschick gave examples that set some limits on how far such a
theorem might be generalized. Firstly, they show that there are elements $a \in \Mod(S)$ which
are not pseudo-Anosov, but nevertheless for which there are homogeneous
quasimorphisms which do not vanish on $a$. Secondly and conversely, they show that
there are elements $b \in \Mod(S)$
of infinite order, for which no power of $b$ is conjugate to its inverse, but for
which $h(b)=0$ for any homogeneous quasimorphism $h$ on $\Mod(S)$.
These examples indicate the subtlety of the function $\scl$ on $\Mod(S)$, and
suggest that a multiplicity of approaches are necessary to appreciate its full
richness.

See \cite{Endo_Kotschick:invent} and \cite{Endo_Kotschick} respectively for
details.

\medskip

Our main theorem about amalgamated free products is:

\begin{amalgam_thm}[Amalgamation Theorem]
Let $G=A*_C B$.
Let $w$ be a word on $A,B$ which is
reduced and cyclically reduced
such that $|w| >1$, and let $\overline{w}$ denote the corresponding element in $G$.
Suppose that the double coset $C\overline{w}C$,
does not contain the element corresponding to any cyclic conjugate of $w^{-1}$.
Then there exists a homogeneous quasimorphism $h$ on $G$
such that $h(\overline{w})=1$, and the defect of $h$ is $\le 312$.
\end{amalgam_thm}

One can make a similar statement about HNN extensions.

Here a word $w$ is said to be {\em cyclically reduced} if 
$w$ is not equal to $v_1 v_2 v_1^{-1} (v_1 \not= \emptyset)$
as a word. The double coset condition is at first glance somewhat odd. We examine a
particular example (of essential slopes in knot complements in $3$-manifolds)
in some detail and show why it is natural and unavoidable.

\subsection{Acknowledgments}

We would like to thank Nathan Dunfield, Benson Farb, \'Etienne Ghys, Daniel Groves, 
Ursula Hamenst\"adt, Dieter Kotschick, Jason Manning, Shigenori Matsumoto, Igor Mineyev,
Shigeyuki Morita, Leonid Polterovich
and Bill Thurston. The second author appreciates the hospitality of 
the department of Mathematics at Caltech, where he visited
the first author, and MPI in Bonn.

While writing this paper, the first author was partially supported by a Sloan Research
Fellowship, and NSF grant DMS-0405491.

\section{Background Material}\label{background_section}

\subsection{$\delta$-hyperbolic groups and spaces}\label{hyperbolic_background_subsection}

We assume the reader is familiar with basic elements of the theory
of $\delta$-hyperbolic groups and spaces: quasi-isometries, quasigeodesics, word-hyperbolicity,
$\delta$-thin triangles, translation length, etc. For a reference, see \cite{Gromov}.

Let $G$ be a $\delta$-hyperbolic group with a fixed generating set $S$.
We denote the translation length of an element $g \in G$ on the
Cayley graph $\Gamma_S(G)$ by $\tau(g)$. This length
is defined by
$$\tau(g) = \lim_{n \to \infty} \frac {\text{dist}(\id,g^n)} n$$
where distance is measured in the usual way in $\Gamma_S(G)$. Observe that
$\tau$ is a class function.

If $g$ is torsion, $\tau(g)=0$. Otherwise, recall that there is a constant
$N(|S|,\delta)$ and an integer $n < N$ such that a power $g^n$ of $g$
leaves invariant a bi-infinite geodesic {\em axis} and
acts on this axis as translation through a distance $n\tau(g)$.
Since this distance is necessarily a positive integer, this implies that there is
a positive constant $B(|S|,\delta)$ such that $\tau(g)\ge B$ for every 
non-torsion element $g$. See e.g. Theorem 5.1 of \cite{epstein-fujiwara}.

\subsection{Stable commutator length}

\begin{definition}
Let $G$ be a group, and $a \in G$. The {\em commutator length} of $a$, denoted
$\cl(a)$, is the minimum number of commutators whose product is equal to $a$.
If $a$ is not in $[G,G]$, we set $\cl(a)=\infty$. The {\em stable commutator
length}, denoted $\scl(a)$, is the lim inf
$$\scl(a) = \liminf_{n \to \infty} \frac {\cl(a^n)} n$$
Notice $\scl(a) = \infty$ if and only if the order of $[a]$ in $H_1(G;\Z)$ is infinite.
\end{definition}

The functions $\cl$ and $\scl$, like $\tau$, are class functions. Observe
that $\cl$ and $\scl$ are {\em monotone} under homomorphisms. That is, if $\phi:G \to H$
is a homomorphism, then
$$\scl(\phi(a)) \le \scl(a)$$
for all $a \in G$, and similarly for $\cl$.

\begin{example}[Mirror Condition]\label{inverse_conjugate_cancels}
Let $G$ be a group, and suppose there are elements $a,b \in G$ and
integers $n \ne m$ such that
$ba^nb^{-1} = a^m$. Then
$$n\cdot \scl(a) = \scl(a^n) = \scl(ba^nb^{-1}) = \scl(a^m) = m\cdot \scl(a)$$
and therefore $\scl(a)=0$. (Note that the defining property of $a$ means
that $a^{n-m} \in [G,G]$, so $\scl(a) \ne \infty$).

Such elements can be found in {\em Baumslag-Solitar groups}. We compute
$$|n|\cdot \tau(a) = \tau(a^n) = \tau(ba^nb^{-1}) = \tau(a^m) = |m|\cdot\tau(a)$$
so if $\tau(a)$ is not zero (which can only happen for a torsion element when $G$ is hyperbolic), 
this is only possible if $n = \pm m$. The group
$$\langle a,b \; | \; ba^{-1}b^{-1} = a , b^2 = \id \rangle$$ may be thought of as the
(orbifold) fundamental group of the {\em interval with mirrored endpoints},
which is sometimes just called a {\em mirror interval}. In a hyperbolic manifold,
an element of $\pi_1$ conjugate to its inverse is represented by a geodesic
which ``double covers'' a geodesic segment with both
endpoints on an orbifold stratum of order $2$ (i.e. it double covers a mirror interval). 
Thus we will sometimes say that
an element $a \in G$ for which there is no $b \in G$ and no integer $n\ne 0$
for which $ba^n b^{-1} = a^{-n}$ satisfies the {\em no mirror condition}.

Analogues of this condition will occur in the hypotheses of all our main
theorems throughout this paper.
\end{example}

\begin{example}\label{hyperbolic_surgery_example}
Let $K$ be a knot complement in $S^3$ with genus $g$, and let $M_{p,q}$ be the
result of $(p,q)$ Dehn filling on $K$. Then if $a \in \pi_1(M_{p,q})$ 
represents the core of the added solid torus, $\scl(a) \le g/p$
(in fact, by taking covers and tubing boundary components together, 
one can actually obtain $\scl(a) \le (g - 1/2)/p$; see \cite{Calegari_length}).
If $K$ is not a satellite knot or a torus knot, for all but finitely
many pairs $(p,q)$ the resulting manifold $M_{p,q}$ is hyperbolic, and
its fundamental group is word-hyperbolic. (See e.g. \cite{Thurston_notes}).
\end{example}

\subsection{Quasimorphisms}\label{qm_subsection}

\begin{definition}
Let $G$ be a group. A {\em quasimorphism} is a function
$$\phi:G \to \R$$
for which there is a least constant $D(\phi)\ge 0$ called the {\em defect},
such that
$$|\phi(a) + \phi(b) - \phi(ab)| \le D(\phi)$$
for all $a,b \in G$. In words, a quasimorphism fails to be linear
by a bounded amount.

A quasimorphism is {\em homogeneous} if $\phi(a^n) = n\phi(a)$
for all integers $n$ and all $a \in G$.
\end{definition}

If $\phi$ is a quasimorphism on $G$, then one can obtain a 
homogeneous quasimorphism $\overline{\phi}$ by the formula
$$\overline{\phi} (a)=\lim_{n \to \infty} \frac{\phi(a^n)}{n}.$$
Note that the defining property of a quasimorphism (that it is ``almost linear'')
implies that the limit exists. A homogeneous quasimorphism is a class function.

The defect $D(\overline{\phi})$ is related to $D(\phi)$ by
$$D(\overline{\phi}) \le 2\cdot D(\phi)$$
see e.g. \cite{Calegari_scl_monograph}, Corollary~2.59.
Homogeneous quasimorphisms are often easier to work with than ordinary quasimorphisms,
but ordinary quasimorphisms are easier to construct. We use this averaging procedure to
move back and forth between the two concepts.

We denote the vector space of all homogeneous quasimorphisms on $G$ by $Q(G)$. Observe
that for any commutator $[a,b] \in G$ and any $\phi \in Q(G)$ we have
$$|\phi([a,b])| \le D(\phi)$$
It turns out that there is an equality
$$\sup_{a,b} |\phi([a,b])| = D(\phi)$$
see \cite{Bavard}, Lemma 3.6.

Quasimorphisms and stable commutator length are related by Bavard's Duality Theorem
(c.f. \cite{Bavard}):

\begin{theorem}[Bavard's Duality Theorem]\label{Bavard_duality_theorem}
Let $G$ be a group. Then for any $a \in [G,G]$, we have an equality
$$\scl(a) = \frac 1 2 \sup_{\phi \in Q(G)} \frac {|\phi(a)|} {D(\phi)}$$
\end{theorem}

Note that one must take the supremum over $\phi \in Q(G)$ with $\phi(a)\ne 0$ 
(and therefore $D(\phi)>0$) for this to make sense; also, if $Q(G) = H^1(G)$ then
$\scl(a) = 0$ for $a \in [G,G]$. Note further that the theorem makes
sense and is true for $a$ satisfying $a^n \in [G,G]$ for some positive $n$.

Bavard's theorem depends on the Hahn-Banach theorem and $L^1-L^\infty$ duality. 
Note that the Hahn-Banach theorem is equivalent to the Axiom of Choice. In particular,
the quasimorphisms promised by Bavard's theorem are typically {\em not constructible}.
Therefore we take explicit note in the sequel of when our arguments make use of Bavard's theorem,
and when they do not.

\begin{example}
For non-elementary
hyperbolic groups $G$, the space $Q(G)$ has an uncountable dimension (\cite{epstein-fujiwara}).
But for certain groups, one can completely understand $Q(G)$. For example, let
$G$ denote the universal central extension of the group of all 
orientation-preserving homeomorphisms of $S^1$. Then $Q(G)= \R$, generated
by Poincar\'e's rotation number. In particular, every non-negative real number is equal to
$\scl$ of some conjugacy class in this group. Similarly, let $\widehat{T}$ denote the
universal central extension of Thompson's group $T$ of dyadic piecewise linear 
homeomorphisms of $S^1$. Then $\widehat{T}$ is a finitely presented group which realizes
every non-negative rational number as $\scl$. See \cite{Ghys_Sergiescu} for more
about the group $T$ and its bounded cohomology.
\end{example}

\subsection{Mineyev's geodesic flow space}

To understand stable commutator length in word-hyperbolic groups, one
needs to control the geometry of maps of surfaces into $\delta$-hyperbolic
spaces. Naively, following the usual practice in hyperbolic manifolds,
one triangulates a surface and ``straightens'' the simplices, and then
appeals to the Gauss-Bonnet theorem to control area and therefore diameter
in terms of injectivity radius. In a $\delta$-hyperbolic space, the
straightening must be done in a careful way. There are at least three
technical approaches to this straightening:

\begin{enumerate}
\item{Gromov's geodesic flow space (\cite{Gromov}, chapter 8)}
\item{Mineyev-Monod-Shalom's homological ideal bicombing (\cite{Mineyev_Monod_Shalom})}
\item{Mineyev's geodesic flow space (\cite{Mineyev})}
\end{enumerate}

The approach in \cite{Gromov} is not entirely fleshed out, and the
``geodesic flow'' is really a quasigeodesic flow; therefore for our applications,
Mineyev's flow space is best suited.

\medskip

Mineyev constructs from a hyperbolic complex $X$ (for instance, a Cayley graph
for a word-hyperbolic group) a geodesic flow space
$\F(X)$, which consists of a union
of parameterized lines joining ordered  pairs of distinct points in
the ideal boundary $\partial X$. The space $\F(X)$ admits a number of
metrics and pseudo-metrics, of which the pseudo-metric $d^\times$ is most important
to us.

The following summarizes some of the main properties of $\F(X)$ which we use.

\begin{theorem}[Mineyev]\label{Mineyev_theorem}
Let $X,d_X$ be a $\delta$-hyperbolic complex with valence $\le n$.
Then there exists a pseudo-metric space $\F(X),d^\times$ called the {\em flow space of $X$}
with the following properties:
\begin{enumerate}
\item{$\F(X)$ is homeomorphic to $(\partial X \times \partial X -\Delta)\times \R$.
The factors $(p,q,\cdot)$ under this homeomorphism are called the {\em flowlines}.}
\item{There is an $\R$-action on  $\F(X)$ (the {\em geodesic flow}) which acts as
an isometric translation on each flowline $(p,q,\cdot)$.}
\item{There is a $\Z/2\Z$ action $x \to x^*$ which anti-commutes with the $\R$ action, 
which satisfies $d^\times(x,x^*)=0$, and which interchanges the 
flowlines $(p,q,\cdot)$ and $(q,p,\cdot)$}
\item{There is a natural action of $\isom(X)$ on
$\F(X)$ by isometries. If $g \in \isom(X)$ is hyperbolic with fixed points $p^\pm$ in $\partial X$
then $g$ fixes the flowline $(p^-,p^+,\cdot)$ of $\F(X)$ 
and acts on it as a translation by a distance which we denote $\tau(g)$. This action of $\isom(X)$
commutes with the $\R$ and $\Z/2\Z$ actions.}
\item{There are constants $M\ge 0$ and $0 \le \lambda < 1$ such that
for all triples $a,b,c \in \partial X$, there is a natural isometric parameterization of
the flowlines $(a,c,\cdot),(b,c,\cdot)$ for which there is exponential convergence
$$d^\times((a,c,t),(b,c,t)) \le M \lambda^t$$}
\item{If $X$ admits a cocompact isometric action, and $G \subset \isom(X)$ is torsion
free, then there is a $G$-equivariant $(K,\epsilon)$ quasi-isometry
between $\F(X),d^\times$ and $X,d_X$.}
\end{enumerate}
Moreover, all constants as above depend only on $\delta$ and $n$.
\end{theorem}

Note that although $\F(X)$ is homeomorphic to $(\partial X \times \partial X - \Delta)\times \R$,
this topology is not induced by the pseudo-metric $d^\times$, since this pseudo-metric fails to
separate pairs of points interchanged by the $\Z/2\Z$ action.

This is a conflation of several results in \cite{Mineyev}. The pseudo-metric $d^\times$
is defined in Section 3.2 and Section 8.6 on a slightly larger space which Mineyev
calls the {\em symmetric join}. The flow space, defined in Section 13, is a natural
subset of this. The basic properties of the $\R,\Z/2\Z$ and $\isom(X)$ action
are proved in Section 2. The remaining properties are subsets of Theorem 44
(page 459) and Theorem 57 (page 468). There is another natural metric
$d_*$ on $\F(X)$ which induces the topology on $\F(X)$, and for which the $\R$ action
is by bi-Lipschitz homeomorphisms.

For our applications, the key points are that the action of hyperbolic isometries
on the flowlines is by translations, and the exponential convergence of flowlines
with common endpoints at infinity.

\begin{remark}
Since the flowlines $(p,q,\cdot)$ and $(q,p,\cdot)$ are distance $0$ apart
in $d^\times$, we may think of them as different parameterizations of the same
lines in a suitable quotient on which $d^\times$ is a metric (and not just a
pseudo-metric).
\end{remark}

\begin{remark}\label{synchronous_parameterizations}
In bullet (5), the synchronous parameterizations of $(a,c,t)$ and $(b,c,t)$ 
for which exponential convergence holds
are precisely those for which $(a,c,0)$ is the point on the flowline
$(a,c,\cdot)$ closest to $b$, and $(b,c,0)$ is the point on the
flowline $(b,c,\cdot)$ closest to $a$ (interpreted in terms of suitable
horofunctions). In particular, for every triple $a,b,c$ of distinct points
in $\F(X)$, if $\Delta$ is a triangle obtained from the union of three flowlines joining
these points in pairs, the edges of $\Delta$ are exponentially close to each other away
from a compact subset of $\Delta$ of uniformly bounded diameter (i.e. the
diameter of the ``thick part'' of $\Delta$ is bounded independently of the choice of $a,b,c$).
\end{remark}

\section{The Gap Theorem, first version}

The Gap Theorem, to be proved below in its first version, says that in a word-hyperbolic group $G$,
if $a$ is an element which satisfies the no mirror condition (from Example~\ref{inverse_conjugate_cancels}) 
and has a sufficiently long translation length, then the stable commutator length
of $a$ can be {\em uniformly} bounded from below. Example~\ref{inverse_conjugate_cancels} and
Example~\ref{hyperbolic_surgery_example} together point to the necessity of
both hypotheses.

\begin{gap_thm}[Gap Theorem, weak version]
Let $G$ be a word-hyperbolic group that is $\delta$-hyperbolic with respect to
a symmetric generating set $S$ with $|S|$ generators. 
Then there is a constant $C(\delta,|S|)>0$
such that for every $a \in G$, either $\scl(a) \ge C$ or else there is
some positive integer $n$ and some $b \in G$ such that $ba^{-n}b^{-1} = a^n$.
\end{gap_thm}
\begin{remark}
By Bavard's theorem, the uniform lower bound on $\scl(a)$ is equivalent to
the existence of a homogeneous quasimorphism $\phi$ with $\phi(a)=1$ for
which there is a uniform upper bound on the defect $D(\phi)$ depending only
on $\delta$ and $|S|$.
\end{remark}

\begin{remark}
If $ba^{-n}b^{-1} = a^n$ for some $b$ and positive integer $n$, then $b^2$ and $a^n$
commute, and since $G$ is hyperbolic, suitable powers of $a$ and $b$ generate
an infinite dihedral group.

It follows that if $G$ is torsion free, no such element $b$ can exist. Conversely,
note that every torsion element $a$ satisfies $a^n = a^{-n} = \id$ for some positive $n$.
It follows that the no mirrors condition is vacuously satisfied in a torsion-free
hyperbolic group.
\end{remark}

The following proof uses properties of
Mineyev's geodesic flow space. We prove a stronger theorem in \S \ref{word_hyperbolic_section}, 
by directly constructing sufficiently many quasimorphisms. The
second proof is logically superior to the first, since the construction is
direct, and moreover the existence of these quasimorphisms does not depend on
Bavard's Theorem and the Axiom of Choice. However, the first proof is more ``geometric''.

\begin{proof}
Let $\Gamma_S(G)$ denote the Cayley graph of $G$ with respect to the generating set $S$.
For any element $a \in G$ recall that $\tau(a)$ denotes the translation length of $a$.
As observed in \S~\ref{hyperbolic_background_subsection}, every element $a \in G$ is either 
torsion, or has a power $a^n$ where $n \le C_2(|S|,\delta)$ that fixes an axis $l_a$, and
therefore has $\tau(a) \ge C_1(|S|,\delta)$. For the remainder of the proof we replace $a$
by a suitable power $a^n$, and assume that $a$ fixes an axis $l_a$.

For convenience, we let $K$ be a $2$-complex coming from a finite presentation
for $G$ with generating set $S$. Observe that $\Gamma_S(G)$ is the $1$-skeleton of
the universal cover $\til{K}$.
Suppose $\scl(a)$ is very small, so that there is an expression
$$a^n = \prod_{i=1}^m [b_i,c_i]$$ 
for some $b_i,c_i \in G$ such that $m/n$ is small.
Let $\Sigma$ be a surface of genus $m$, and $\sigma:\Sigma \to K$ a simplicial map
which takes the boundary to the immersed circle 
$\gamma = l_a/a^n$. We would like to choose a representative surface whose
area and geometry can be controlled by $\chi(\Sigma)$ and $\length(\gamma)$. Since we need
to control constants, this must be done carefully.

\medskip

We appeal to Mineyev's Theorem~\ref{Mineyev_theorem}. Using this theorem, one proceeds as follows.
One picks a $1$-vertex triangulation of $\Sigma$, and spins the vertex along
$\partial \Sigma$ thereby producing an ideal triangulation whose edges
can be realized as a subset in the quotient space $\F(\Gamma_S(G))/G$. 
This is covered by a $\pi_1(\Sigma)$ equivariant map
from the $1$-skeleton of $\til{\Sigma}$ to $\F(\Gamma_S(G))$,
where $\pi_1(\Sigma)$ acts on $\F(\Gamma_S(G))$ by its image in $G$.
Each ideal triangle in $\til{\Sigma}$ corresponds to a
union of three flowlines in $\F(\Gamma_S(G))$. 

In $\F$, as pointed out in Remark~\ref{synchronous_parameterizations},
each triangle consists of a thick {\em core} of
diameter at most $C_3$ (depending only on $\delta$ and $|S|$) 
together with thin regions consisting of pairs of geodesic
rays whose distance converges exponentially fast to zero. The image of
$\Sigma$ consists of $2|\chi(\Sigma)| = 4m-2$ triangles. Because of the
exponential convergence of geodesics in $\F$,
we can give $\Sigma$ a hyperbolic metric of constant curvature
$-\kappa$ in such a way that the triangles are totally geodesic, and the
map to $\F(\Gamma_S(G))$ is distance decreasing on the $1$-skeleton, 
where $\kappa>0$ depends only on
$\delta$ and $|S|$. It follows from the Gauss-Bonnet theorem that
away from a thick part consisting of at most $(4m-2)$ regions whose
diameters sum to at most $(4m-2)C_3$ for some constant $C_3$, 
the thickness of $\Sigma$ is bounded by a constant $C_4$,
where $C_3$ and $C_4$ depend only on $|S|$ and $\delta$. In fact, by
choosing $C_3$ sufficiently large, we may assume $C_4$ is as small as we like,
a fact which we will not use.

After composing with a
quasi-isometry $$\phi:\F(X) \to \til{K}$$ and filling in the map on
triangles, we get an induced map
$$\til{\sigma}_1:\til{\Sigma} \to \til{K}$$
It doesn't really matter how the map is filled in on triangles, since we
are ultimately only interested in the distances between points contained in the
image of $\partial \til{\Sigma}$.
Again,  by Theorem~\ref{Mineyev_theorem}, we can assume that 
the quasi-isometry constants of $\phi$ depend only
on $\delta$ and $|S|$.

Under this quasi-isometry, the constants $C_3$ and $C_4$ for the diameter of
the core and the thickness of the complementary region
must be replaced by analogous constants $C_3'$ and $C_4'$ which still depend only
on $\delta$ and $|S|$.
The image of every boundary component is a quasigeodesic which
is within a bounded distance of some translate of the axis $l_a$;
it follows that for a suitable choice of $\phi$,
without changing the constants involved, we can 
assume that $\til{\sigma}_1$ takes $\partial \til{\Sigma}$
in a $\pi_1(\Sigma)$-equivariant way to a union of translates of $l_a$.

By abuse of notation, we denote the image of this map as $\til{\Sigma}_1$.
Let $\alpha$ be a fundamental domain for the action of
$a^n$ on $l_a$. By our estimates, away from a subset of $\alpha$ of length at most
$(4m-2)C_3'$, every point $p \in \alpha$ can be joined by an arc $\til{\beta}_p$
in the image of $\til{\Sigma}_1$,
such that the endpoints of $\til{\beta}_p$ are a distance at most $C_4'$
apart in $\til{K}$, and lie on distinct components $l_a$ and $l_a^i$ of 
$\til{\sigma}_1(\partial \til{\Sigma}_1)$. Pulling
back by $\til{\sigma}_1$ and projecting to $\Sigma_1$, we obtain a homotopically essential
proper arc $\beta_p \in \Sigma_1$. If $p_1,p_2$ are at least distance $2C_4'$ apart
in $\alpha$, we see that
the pulled-back arcs $\beta_{p_1},\beta_{p_2}$ can be isotoped to be
disjoint in $\Sigma_1$.
A surface of genus $m$ with one boundary component contains at most $6m-3$ disjoint
nonparallel proper essential arcs, and each arc has $2$ endpoints.
It follows that there are at most $12m-6$ components $l_a^i$ which are joined by
arcs $\til{\beta}_p$ with nearby endpoints to points in $\alpha$.

For a fixed component $l_a^i$, the set of points where $l_a^i$ is close to $l_a$
is (coarsely) connected, by convexity of quasigeodesics in $\delta$-hyperbolic spaces.
Consequently there is a connected arc $\sigma' \subset \alpha$
satisfying
$$\length(\sigma') \ge \frac {\length(\alpha) - (4m-2)(C_3'+ C_4')} {12m-6}$$
which cobounds a strip $R'$ of $\til{\Sigma}_1$ of width $\le C_4'$ with some
fixed $l_a^i$. 

By convexity of quasigeodesics, there is a connected subinterval $\sigma \subset \sigma'$
satisfying
$$\length(\sigma) \ge \length(\sigma') - 2C_4'$$
such that $\sigma$ cobounds a strip $R$ of $\til{\Sigma}_1$ of width
$C_5$ depending only on $\delta$ with some fixed $l_a^i$. The strip $R$ is much longer
than it is wide, and it makes sense to say that a choice of orientations on the
$l_a$ and $l_a^i$ agree or disagree (see Fig. \ref{agree_figure}) :

\begin{figure}[ht]
\center{\scalebox{.5}{\includegraphics{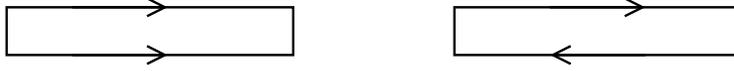}}}
\caption{The orientations on $l_a$ and $l_a^i$ might agree or disagree along $R$}\label{agree_figure}
\end{figure}

\medskip

Now, since $\Sigma_1$ is oriented, the induced orientations
on $l_a$ and $l_a^i$ {\em disagree}. Therefore
$l_a^i$ is a nontrivial translate of the axis $l_a$. That is, there is some $c$
with $l_a^i = c(l_a)$. Then
$e=ca^{-1}c^{-1}$ stabilizes $l_a^i$, and moves points in (roughly)
the same direction as $a$.

Let $p$ be the midpoint of $\sigma$. Suppose that $\length(\gamma) = \tau(a)$
is big compared to the constant $C_5$ (e.g. $\length(\gamma) \ge 2C_5$ will do).
We consider the translates
$e^{-w}a^w(p)$ for
$$|w| \le \frac {\length(\sigma)} {3 \cdot \length(\gamma)} \le \frac n {36m} + 1$$
where the last inequality follows from our assumptions about
$\length(\gamma)$.

Then by our estimate on the width of the strip $R$, we have
$$d(p,e^{-w}a^w(p)) \le 4C_5$$

There are at most $|S|^{4C_5}$ elements in the group which move the point $p$
a distance $\le 4C_5$, and therefore if
$$\frac n {36m} + 1> |S|^{4C_5}$$
then by the pigeonhole principle there are $w_1 \ne w_2$ such that
$$e^{-w_1}a^{w_1} = e^{-w_2}a^{w_2}$$
and therefore
$$ca^{w_1-w_2}c^{-1} = e^{w_2-w_1} = a^{w_2-w_1}$$
In words: if $\scl(a)$ is too small,
a nontrivial power of $a$ is conjugate to its inverse.

It follows contrapositively that when $\tau(a) \ge 2C_5$ and no power of $a$ is conjugate to
its inverse, we have $\scl(a) \ge \frac 1 {36} |S|^{-4C_5}$.
If $a$ is not torsion, $\tau(a) \ge C_1$. So for any $a$ with no power conjugate
to its inverse, we have an estimate 
$$\scl(a^{2C_5/C_1}) \ge 
\frac 1 {36} |S|^{-4C_5}$$
And therefore the constant $$C(\delta,|S|) = \frac 1 {36} |S|^{-4C_5} \frac {C_1} {2C_5}$$
satisfies the claim of the theorem.
\end{proof}

\begin{remark}
It follows from the proof that if $a$ does not satisfy the no mirror condition, so that
there is a positive $n$ with $a^n = ca^{-n}c^{-1}$ for some $c$, 
then $n$ may be chosen to be less than some constant depending only on $\delta$ and $|S|$.
\end{remark}

It is interesting to ask the following:

\begin{question}\label{gap_question}
Let $G$ be a group, and let $H$ be the set of elements in $G$
which satisfy the no mirror condition; i.e. for which $a^n \ne ba^{-n}b^{-1}$ whenever
$b\in G$ and $n \ne 0$.
Under what circumstances is $\inf_{a \in H} \scl(a)$ positive?
\end{question}

We say that a group $G$ has the {\em stable commutator gap property} if it has the property
described in Question~\ref{gap_question}. 
Theorem A shows that hyperbolic groups have the stable commutator gap property.

\begin{example}\label{lower_bound_example}
In \cite{Comerford_Edmunds}, Comerford and Edmunds show that in a free group $F$,
every nonzero element $a$ satisfies $\scl(a) \ge 1/2$ (this lower bound essentially follows
from an earlier result of Duncan and Howie \cite{Duncan_Howie} and improves an earlier
estimate $\scl(a) \ge 1/6$ obtained by Culler in \cite{Culler}).
Let $G$ be residually free; i.e. suppose for all nonzero $a \in G$ there is
a homomorphism $\phi_a:G \to F$ to a free group, for which $\phi_a(a)$ is nonzero. Since $\scl$
is nonincreasing under homomorphisms, $G$ satisfies the stable commutator gap property.
\end{example}

\section{Counting functions}\label{counting_function_section}

Here we review so-called {\em counting functions}, which generalize a construction 
introduced by Brooks \cite{Brooks}
to construct quasimorphisms on free groups. See  
\cite{epstein-fujiwara}, \cite{fujiwara:plms} or \cite{Fujiwara_Bestvina} for more details.

Suppose $G$ is a group with a fixed symmetric generating set $S$, and $\Gamma:=\Gamma_S(G)$ is
its Cayley graph. Let $w$ be a word in the generating set.
Let $\alpha$ be a (directed) path in $\Gamma$,
and $|\alpha|$ its length.
Define $|\alpha|_w$ to be the maximal
number of times that $w$ can be seen as an (oriented) subword of $\alpha$ without overlapping.

\begin{example}
$|xyxyx|_{xy} = 2$. $|xyxyx|_{xyx} = 1$. $|xxyxy|_{yx} = 1$.
\end{example}

Let $0<W<|w|$ be a constant.
For $x,y \in \Gamma$, 
define
$$c_{w,W}(x,y):=d(x,y)
- \inf_{\alpha}(|\alpha|-W|\alpha|_w),$$
where $\alpha$ ranges over all the paths from $x$ to $y$.
If the infimum is attained by $\alpha$, we 
say $\alpha$ is a {\em realizing path} for $c_{w,W}$ from $x$ to $y$.
If $\gamma$ is a geodesic from $x$ to $y$, then define
$c_{w,W}(\gamma)=c_{w,W}(x,y)$.

Fix a point $x \in \Gamma$.
Define for $a \in G$
$$c_{w,W}(a):=c_{w,W}(x,ax).$$
$c_{w,W}$ is called the {\em counting function}
for the pair $(w,W)$.
Let $w^{-1}$ denote the inverse word of $w$.
We define
$$h_{w,W}:=c_{w,W}-c_{w^{-1},W}.$$
In \cite{epstein-fujiwara}, the normalization $W=1$ is used. 
This is an appropriate choice of
constant when $w^*:=\cdots wwww \cdots$ is a bi-infinite geodesic,
which is the case in our applications throughout this paper.

\medskip

More generally, when $G$ acts on a graph $\Gamma$ (not necessarily properly), 
we modify the definition of counting functions as follows.
Let $w$ be a path in $\Gamma$ and 
call $aw$ for $a \in G$ a {\em copy of $w$} or (interchangeably) a {\em translate of $w$}.
For a path $\alpha$ in $\Gamma$, 
define $|\alpha|_w$ to be the maximal
number of disjoint oriented copies of $w$ which can be obtained as subpaths of $\alpha$.
All other definitions are as above.
These modified counting functions will be important in the sequel when we 
discuss the action of the mapping class
group $\Mod(S)$ on the complex of curves $\CC(S)$, for some surface $S$.

\begin{prop}[Lemma 3.3 \cite{fujiwara:plms}]\label{quasigeodesic.realizing}
If $\alpha$ is a realizing path
for $c_{w,W}$, then it is a $(K,\epsilon)$-quasigeodesic, where
$$K=\frac {|w|} {|w|-W} \; , \; \epsilon = \frac {2W|w|} {|w|-W}$$
\end{prop}

It is known that in a $\delta$-hyperbolic space, 
any two $(K,\epsilon)$-quasigeodesics which 
have same end points stay in the $L(K,\epsilon,\delta)$-neighborhood 
of each other for some universal constant $L$ (see \cite{Gromov}).

Let $L=L(|w|/(|w|-W),2W|w|/(|w|-W),\delta)$.
Let $\alpha$ be a geodesic from $x$ to $y$.
From Proposition~\ref{quasigeodesic.realizing} we deduce
that a realizing path for $\alpha$ must be contained in the $L$-neighborhood of $\alpha$.
Consequently, if the $L$-neighborhood of $\alpha$
does not contain a copy of $w$, then
$c_{w,W}(\alpha)=0$.

\begin{remark}\label{geodesic_remark}
We will use this fact later in our argument to avoid ``reverse counting''. 
Roughly speaking, let $w$ be a word such that $w^n$ is a geodesic.
Then, for $n>0$, 
$$c_{w,W}(w^n) \ge Wn$$ 
because $|w^n|_w=n$.

Suppose the $L$-neighborhood of $w^n$ does not
contain a copy of $w^{-1}$. Here we are thinking of the
$L$-neighborhood of $w^n$, for large $n$, like a long
narrow tube whose core has a definite orientation, agreeing
with the orientation on $w$. By ``a copy of $w^{-1}$'' here we
mean a copy of $w$ whose orientation {\em disagrees} with that
of the core of the tube (compare with Fig. \ref{agree_figure} and the
accompanying discussion).
Then it follows that $c_{w^{-1},W}(w^n)=0$
because for a realizing path $\alpha$
for $c_{w^{-1},W}$ at $w^n$ we must have
$|\alpha|_{w^{-1}}=0$.
We thus obtain for all $n>0$ an inequality
$h_{w,W}(w^n) \ge nW$.
\end{remark}

Let $D(h)$ be the defect of $h$. Then we have the following inequality:

\begin{prop}[Prop 3.10 \cite{fujiwara:plms}]\label{defect_estimate}
$$D(h_{w,W}) \le 12 L +6W+48 \delta.$$
\end{prop}

\begin{remark}\label{subsequent_remark}
Note that the defect only depends on $|w|,W,\delta$.
If we take $W=1$, then $L$ depends only on $\delta$ if 
$|w| \ge 2$. 
In particular, the upper bound in 
Proposition~\ref{defect_estimate} depends {\em only} on $\delta$.
We will take $W=1$ in all the applications in this paper.
\end{remark}

\section{The Gap Theorem, improved version}\label{word_hyperbolic_section}

We are now in a position to state and prove the improved version of the Gap Theorem,
using counting quasimorphisms.
The following theorem improves Theorem A in at least two ways: it does not
use Bavard's theorem or Mineyev's flow space, and the constants in questions can
be effectively estimated.

\begin{improved_gap_thm}[Gap Theorem, strong version]
Let $G$ be a word-hyperbolic group that is $\delta$-hyperbolic
with respect to a symmetric generating set $S$ with $|S|$ generators.
Let $a \in G$ be a (non-torsion) element.
Assume there is no $n >0$ and no $b \in G$
with 
$ba^{-n}b^{-1}=a^n$.
Then there is a homogeneous quasimorphism $h$ on $G$
such that 
\begin{enumerate}
\item{$h(a) =1 $}
\item{The defect of $h$ is $\le C(\delta,|S|)$.}
\end{enumerate}
Moreover, let $a_i \in G$ be a collection of elements for which
$T = \sup_i \tau(a_i)$ is finite. Suppose that for all integers $n,m \ne 0$
and all elements $b \in G$ and indices $i$, there is an inequality
$$ba^nb^{-1} \ne a_i^m$$
Then there is a homogeneous quasimorphism $h$ on $G$
such that 
\begin{enumerate}
\item{$h(a) =1 $, and $h(a_i)=0$ for all $i$}
\item{The defect of $h$ is 
$\le C'(\delta, |S|, T)$}
\end{enumerate}
\end{improved_gap_thm}

\begin{remark}
Note that in the second part of the theorem, we must state a generalization of
the no mirror condition, involving both $a$ and the elements $a_i$.
\end{remark}

\begin{remark}
Hamenst\"adt (\cite{Ha}) used dynamical methods to directly construct quasimorphisms on
word-hyperbolic groups. Similar ideas were developed by Picaud (\cite{Pic}) in the
special case of surface groups. It is not clear to us how to use these methods to obtain
uniform estimates on defects.
\end{remark}

\begin{proof}
We follow \cite{epstein-fujiwara}, \cite{fujiwara:plms}.
We assume in what follows that $a$ is not torsion. 
Let $N$ be a constant such that $a^N$ stabilizes a bi-infinite geodesic $l_a$.
Note that one may find an $N$ whose size can be bounded in terms of $\delta$ and $|S|$.
We also choose $N$ so that $N \tau(a) \gg L(\delta)$. Note that since $\tau(a)$ can be bounded
from below by a positive constant depending only on $\delta$ and $|S|$, the constant
$N$ still depends only on $\delta$ and $|S|$
(see the proof of Theorem A).

After replacing $a$ by a conjugate if necessary, we may assume that 
$1 \in l_a$.  
Set  $b=a^N$. $l_a$ is an axis for $b$.
We denote the subpath of $l_a$ from $1$ to $b$
by $l_a|b$ and the subpath from $1$ to $b^{-1}$ by $l_a|b^{-1}$.
Note that we think of these as {\em oriented} segments.

Now let $c_b$ be the counting function on $G$ for the 
pair $(l_a|b,1)$ (note that we are setting $W=1$ in the notation
of \S \ref{counting_function_section}).
Then for any $n>0$, we obtain an estimate $c_b(b^n) \ge n$, since 
the subpath of $l_a$ from $1$ to $b^n$ can be tiled by $n$ disjoint translates
of $l_a|b$. At the cost of possibly replacing $N$ by $2N$ if necessary, we may assume that 
$|l_a|b| \ge 2$, so that
$L$ (as in \S \ref{counting_function_section}) depends only on $\delta$.

\medskip
{\noindent \bf Claim 1:}
There is a constant $C_1(\delta,|S|)$ such that
if $N \ge C_1$ then there is no translate of $l_a|b^{-1}$
in the $L$-neighborhood of $l_a$ whose orientation
agrees with that of $l_a$.
\medskip

Note that $l_a|b^{-1}$ is a copy of $l_a|b$ with the opposite
orientation, so we could just as well state Claim 1 as saying
that there is no translate of $l_a|b$ in the $L$-neighborhood of
$l_a$ whose orientation {\em disagrees} with that of $l_a$.
We will prove the claim later; we call the conclusion of the
claim ``no reverse counting'', and prove our theorem under this hypothesis.

Take $N$ to satisfy  $N \ge C_1$ as well in the following.
Note that $N$ can still be chosen satisfying these criteria
with size bounded from above in terms of $\delta$ and $|S|$.

Let $c_{b^{-1}}$ be the counting function 
for the pair $(l_a|b^{-1},1)$.
It follows from Claim 1 that 
$$c_{b^{-1}}(b^n) =0$$ for all $n >0$.
We define $h_b=c_b - c_{b^{-1}}$, then 
obtain 
$$h_{b}(b^n) \ge n$$ 
for any $n >0$.

\medskip

By Proposition~\ref{defect_estimate} and Remark~\ref{subsequent_remark}, there is 
a constant $K(\delta)$, which depends only 
on $\delta$ such that 
$$D(h_b) \le K(\delta).$$

Since $b=a^{N}$, we get
$h_b(a^n) \ge n/N$ for all $n >0$.
By averaging, we may replace $h_b$ by a homogeneous quasimorphism
$h'=\bar h_b$ (see \S \ref{qm_subsection}).
Then $D(h') \le 4 K(\delta)$ and 
$h'(a) \ge 1/N$.
Define $h=  m h'$ for some constant $m$ for which $h(a)=1$.
Then $h$ is a homogeneous quasimorphism, and satisfies
$D(h) \le 4NK(\delta)$.

This proves the first part of the theorem, modulo Claim 1.
\medskip

The second part follows by almost the same argument. 
First of all, we may assume that $a_i$ has infinite order (equivalently, $\tr(a_i)>0$)
since every homogeneous quasimorphism will already vanish on torsion.
Since we are looking for a homogeneous quasimorphism
which vanishes on $a_i$, without loss of generality,
we may replace each $a_i$ by a conjugate of a non-trivial power
($\le N(\delta,|S|)$).
Therefore, we may assume
that $a_i$ has a geodesic axis $l_{a_i}$, 
with $1 \in l_{a_i}$.

Let $v$ be an oriented geodesic path in $\Gamma_S(G)$ from $1$ to $a^N$, and
let $v^{-1}$ be $v$ with the opposite orientation. In other words, $v = l_a|_{a^N}$.

A relative version of Claim 1 proves the second part of the theorem. We give
two different proofs, which give different constants. Recall our notation $b = a^N$.

\medskip
{\noindent \bf Claim 2:}
Suppose that $N$ satisfies $$N \ge C_1|S|^{\tau(a)}\frac {\tau(a)} {\tau(a_i) +1}$$ 
for suitable $C_1$ depending only on $\delta$ and $|S|$.
Then there is no copy of $l_a|b$ or $l_a|b^{-1}$ in the 
$L$-neighborhood of $l_{a_i}$.

\medskip
{\noindent \bf Claim 2${}^\prime$:}
Suppose that $N$ satisfies $$N \ge C_1|S|^T$$
for suitable $C_1$ depending only on $\delta$ and $|S|$.
Then there is no copy of $l_a|b$ or $l_a|b^{-1}$ in the
$L$-neighborhood of $l_{a_i}$.
\medskip

It follows from Claim 2 or 2${}^\prime$ that $h$ as constructed above satisfies
$h(a_i)=0$ for all $i$.

\medskip

{\noindent \it Proof of Claim 1.}
In fact, the proof of the claim follows by the same argument as the end of the
proof of Theorem A. A copy of $l_a|b^{-1}$ in the $L$-neighborhood of $l_a$ is
contained in an axis of some element $e$ which is conjugate to $a^{-1}$; i.e.
$e = c a^{-1} c^{-1}$ for some $c$. 

Let $v$ denote such a copy of $l_a|b^{-1}$, and let $v_0$
represent its initial point.
The geodesic $l_a$ is invariant under $a^N$ but not necessarily under $a$
itself. Nevertheless, there is a constant $C_2$ depending
only on $\delta$ such that $a^i(v_0)$ is within distance $C_2$ of
$l_a$, for any integer $i$.
Then for all $i < N$ the element $e^{-i}a^i$ satisfies
$$d(v_0,e^{-i}a^i(v_0)) \le 2L+2C_2$$
So if $N>|S|^{2L+2C_2}$ we must have
$e^{-i_1}a^{i_1} = e^{-i_2}a^{i_2}$ for
two distinct indices $i_1,i_2$, 
and therefore the $(i_2-i_1)$-th power of $a$ is conjugate by $c$
to its inverse.
Set $C_1=|S|^{2L+2C_2}$.
\qed

\medskip

{\noindent \it Proof of Claim 2.}
We only sketch the proof of Claim 2 since the details are almost 
identical to those of the proof of Claim 1. If the $L$-neighborhood
of $l_{a_i}$ contains a copy of a arbitrarily long subpath
of $l_a$, we can find a point $v_0$ for which 
$$d(v_0,e^{k_j}a_i^j(v_0))\le 2L + 2C_2 + \tau(a)$$
for some conjugate $e$ of $a$ or $a^{-1}$, and for any $0 \le j \le J$ for
an arbitrarily big (fixed) $J$. It follows that if $N > |S|^{2L+2C_2+\tau(a)}$
then some non-trivial power of $a_i$ is conjugate to some (possibly trivial) power
of $a$, say $a^k$. Moreover, if 
$$N > |S|^{2L+2C_2+\tau(a)} \frac {\tau(a)} {\tau(a_i)+1}$$ 
then $k \not=0$, contrary to hypothesis.

\medskip

{\noindent \it Proof of Claim 2${}^\prime$.}
This argument interchanges the roles of $a$ and $a_i$ in the proof of Claim 2.
Suppose the $L$-neighborhood of $l_{a_i}$
contains a copy, $v$, of $l_a|b$ or $l_a|b^{-1}$.
The segment $v$ is a part of an axis of $e=cac^{-1}$
or $ca^{-1}c^{-1}$. Let $v_0$ be the starting point of $v$.
Then for each $0 \le j \le N$, 
there exists $k_j$ such that 
$$d(v_0,a_i^{k_j} e^j(v_0)) \le 2L+2C_2+\tau(a_i) \le 2L+2C_2+T$$
(we use the same constant $C_2$ for $a_i$).
Therefore, if $N>|S|^{2L+2C_2+T}$, then 
some non-trivial power of $a$ is conjugate to some
(possibly trivial) power of $a_i$, which is 
impossible.

This completes the proof of the theorem.

\end{proof}

\begin{remark}
The second statement can be used
to give a lower bound of the 
stable commutator length ``relative to $\{a_i\}$''. 
That is, if we can write $a$ as a product
$$a = [b_1,c_1]\cdots[b_n,c_n]a_{i_1}^{m_1} \cdots a_{i_k}^{m_k}$$ 
then we say the {\em relative length} of $a$ is $\le n +k/2$.
The infimum of this number is the {\em relative commutator length}, and the liminf of the
relative commutator length of $a^n$ divided by $n$ as $n \to \infty$ is the {\em relative
stable commutator length}.
Using our theorem, the relative commutator 
length  has a lower bound of $1/2D(h)$, where 
$h$ is a homogeneous quasimorphism obtained in the 
second part for $a$ and the $a_i$.
\end{remark}

\begin{remark}
One may reinterpret the second part of Theorem A${}^\prime$ as follows. In any group $G$, let
$B_1(G)$ denote the real vector space of group $1$-boundaries 
(i.e.\ group $1$-cycles that are boundaries of group $2$-chains), and for any chain
$\sum t_i g_i \in B_1(G)$ define $\scl$ by the formula
$$\scl(\sum t_i g_i) = \sup_\phi \frac {|\sum t_i \phi(g_i)|} {2D(\phi)}$$
(compare with Theorem~\ref{Bavard_duality_theorem}). This function is a pseudo-norm, and
agrees with $\scl$ on ordinary elements (see \cite{Calegari_scl_monograph}, \S~2.6). 
It evidently vanishes on the subspace
$H$ spanned by cycles of the form $g - hgh^{-1}$ and $g^n - ng$, and descends to a quotient
pseudo-norm on $B_1^H(G):=B_1(G)/H$. Then the second part of 
Theorem A${}^\prime$ implies that whenever
$G$ is hyperbolic, $\scl$ is a genuine {\em norm} on $B_1^H(G)$.
\end{remark}

We refer to the homogeneous quasimorphisms constructed in the proof of Theorem A${}^\prime$
as {\em counting quasimorphisms}, by contrast with the abstract quasimorphisms
promised by Bavard's theorem.

It is interesting to speculate that one could use Theorem A or A${}^\prime$ as a starting point
to invert the word-hyperbolic Dehn surgery
theory developed by Groves and Manning and independently
by Osin (c.f. \cite{Groves_Manning}, \cite{Osin}).

\begin{corollary}
Let $K$ be a knot in $S^3$ of genus $g$. 
Then for any $\delta>0$ and every integer $n>0$ there
is a constant $C(\delta,n)>0$ such that if $M_{p/q}$ is the result of $p/q$ surgery on $K$,
and $|p| \ge C\cdot g$ then every homomorphism from $\pi_1(M_{p/q})$ to a torsion-free
$n$-generator $\delta$-hyperbolic group is trivial.
\end{corollary}
\begin{proof}
By construction, if $a$ represents the core geodesic of $M_{p/q}$, we can estimate
$$\scl(a) \le g/|p|$$
Moreover, since $K$ is a knot in $S^3$, $\pi_1(M_{p/q})$ is normally generated by $a$.

Let $G$ be $\delta$-hyperbolic and torsion free. Then no element is conjugate to its
inverse, and therefore there is a uniform lower bound on the stable commutator length of
any nontrivial element in $G$, depending only on $\delta$ and the size of
a generating set for $G$. If
$\rho:\pi_1(M_{p/q}) \to G$ is any homomorphism, then
$\scl(\rho(a)) \le \scl(a)$, and the image is normally generated by $\rho(a)$.
The claim follows.
\end{proof}

\begin{remark}
One knows that under a degree $1$ map between hyperbolic manifolds of the same dimension,
volume must go down. It is therefore 
significant in this corollary that the volumes of the manifolds $M_{p/q}$ go
{\em up} as $|p| \to \infty$. 
\end{remark}

\section{First accumulation point}\label{accumulation_section}

In a torsion free word-hyperbolic group $G$, stable commutator length defines a function
$\scl$ from conjugacy classes to $\R$. By Theorem A or A${}^\prime$, the first accumulation point
$\delta_\infty$
for the image of this function is positive. In this section, we obtain {\em universal}
estimates for $\delta_\infty$ which are {\em independent} of $G$.

\begin{accumulation_thm}[Accumulation Theorem]
Let $G$ be a torsion free non-elementary word-hyperbolic group. Then the first accumulation point
$\delta_\infty$ for stable commutator length satisfies
$$\frac 1 {12} \le \delta_\infty \le \frac 1 2$$
\end{accumulation_thm}
\begin{proof}
To obtain the upper bound, observe that $G$ contains a quasi-isometrically embedded
copy of the free group on $2$ generators. A nonabelian free group contains infinitely
many conjugacy classes with $\scl \le 1/2$ (see e.g. \cite{Bavard}).
Since the embedding is quasi-isometric,
the image of infinitely many of these conjugacy classes stay non-conjugate in $G$.
Under any homomorphism, $\scl$ cannot go up, so the upper bound is proved.

An elementary argument gives a lower bound of $1/24$.
We follow the argument in the proof of Theorem A, and
we adopt notation and the setup from that theorem.
In any hyperbolic group, there are only finitely many conjugacy classes of elements whose
translation length is bounded above by any constant. Therefore, we may
assume the translation length of an element $a$ is as long as we like. In particular,
we can assume that there is an axis $l_a$ which is
geodesic and invariant under $a$. As before, let $\alpha$ denote a fundamental
domain for $a^n$. Then $\length(\alpha) = n\tau(a)$.

We suppose, as in the proof of Theorem A, that there is a segment $\sigma \subset l_a$
satisfying
$$\length(\sigma) \ge \frac {\length(\alpha) - (4m-2)(C_3'+C_4')} {12m-6} - 2C_4'$$
which cobounds a strip $R$ of $\til{\Sigma}_1$ of width $\le C_5$ with a translate $c(l_a)$. 
As before, there is $e = ca^{-1}c^{-1}$ which
stabilizes $c(l_a)$ and moves points in almost the same direction as $a$.

If the translation length of $a$ is long enough, then under the
assumption $\scl(a) < 1/24$, the estimate above gives us that
$\length(\sigma)$ is more than twice as big as a fundamental domain for $a$. 
For, if $\scl(a) < 1/24$ then we can choose $n,m$ as above such that $n/m > 24$. Now,
$\length(\alpha) = n \tau(a)$, so
$$\length(\sigma) \ge \frac {n} {12m-6} \tau(a) - O(1) \ge 2\tau(a)$$
providing $\tau(a)$ is sufficiently big.

In fact, we can assume that
$$ \frac {\length(\sigma)} {2\tau(a)} \ge 1+ \epsilon$$
where $\epsilon$ is any number smaller than $1/24 - \scl(a)$.
In particular, we can assume
$$\length(\sigma) - 2\tau(a) = 2C_6$$
where $C_6$ is as big as we like. 
We parameterize $\sigma$ as $\sigma(i)$ where 
$$i \in \left[ \frac {-|\length(\sigma)|} 2, \frac {|\length(\sigma)|} 2 \right]$$
Observe that for all $i \in [-C_6,C_6]$ we have estimates
$$d(e^{-1}a(\sigma(i)),\sigma(i)) \le 4C_5, \; d(ae^{-1}(\sigma(i)),\sigma(i)) \le 4C_5$$
For all $-C_6 \le i \le C_6$ we let $\mu_i \in G$ be such that
$$\mu_i(\sigma(0)) = \sigma(i)$$

It follows that for any $|i| \le C_6$ we have
$$d(\mu_i^{-1}e^{-1}a\mu_i(\sigma(0)),\sigma(0)) \le 4C_5, \; 
d(\mu_i^{-1}ae^{-1}\mu_i(\sigma(0)),\sigma(0)) \le 4C_5$$

If $C_6$ is very big compared to $(|S|^{4C_5})^2$ then by the pigeonhole principle there
are $i_1,\dots, i_n$ where $n \ge C_6/(|S|^{4C_5})^2$ for which
$\mu_{i_j}\mu_{i_k}^{-1}$ commutes with both $e^{-1}a$ and $ae^{-1}$ whenever
$1 \le j,k \le n$. Up to this
point, our argument makes no use of the fact that $G$ is torsion free.

In a torsion free word-hyperbolic group, two nontrivial elements which commute are proportional.
Since $\mu_{i_j}\mu_{i_k}^{-1}$ commutes with both $e^{-1}a$ and $ae^{-1}$, it
follows that $e^{-1}a$ and $ae^{-1}$ themselves are proportional. Since they
are conjugate, they have the same translation length, and 
are either equal or inverse, since $G$ is torsion free. In the first case, $a$ and $e^{-1}$
commute; since they have the same (positive) translation length, they are either equal or
inverse. Since their axes are almost oppositely aligned along $\sigma$, they must
be inverse, so $a=e$. But $e = ca^{-1}c^{-1}$ which is absurd in a torsion free
group. In the second case, $e^{-1}a = ea^{-1}$ so $e^2 = a^2$, and therefore
$e=a$ since $G$ is torsion free, and we get a contradiction in any case. This proves
the estimate $\delta_\infty \ge 1/24$.

\medskip

To get the estimate $\delta_\infty \ge 1/12$ we use Theorem~\ref{Mineyev_theorem}.
We argue as above that if $\scl(a)< 1/12$ and $a$ has sufficiently long translation
length, then $$|\length(\sigma) - \tau(a)| = 2C_6$$ is as big as we like. Let
$p$ be a point at a distance $C_6$ from one of the endpoints of $\sigma$.

Let $p,q$ be the ideal points stabilized by $a$, and $r,s$ the ideal points stabilized
by $b$. Then for suitable parameterizations of the flowlines $(p,q,\cdot),(r,s,\cdot)$,
the points $(p,q,t)$ and $(r,s,t)$ are within distance $Me^{\lambda C_6}$
for $t \in [C_6, \length(\sigma) - C_6]$, where $M$ and $\lambda<1$ are
universal constants (depending only on $\delta$ and $|S|$) but $C_6$ is as big
as we like.

This requires some explanation: by Theorem~\ref{Mineyev_theorem} bullet (5),
for any three ideal points $x,y,z$
there are parameterizations of $(x,z,\cdot)$ and $(y,z,\cdot)$ for
which $$d^\times((x,z,t),(y,z,t)) \le M\lambda^t$$
for suitable $M$ and $0\le \lambda < 1$. In fact, by
Remark~\ref{synchronous_parameterizations}, these parameterizations are exactly those
for which $(x,z,0)$ is the closest point on $(x,z,\cdot)$ to $y$, and $(y,z,0)$ is the
closest point on $(y,z,\cdot)$ to $x$. Now, if $(p,q,\cdot)$ and $(r,s,\cdot)$ are flowlines
which have long sub-segments which are distance $\le C_5$ apart, then if
$m$ is the point on $(p,q,\cdot)$ which is closest to $r$ and $n$ is the point on $(p,q,\cdot)$
which is closest to $s$, then $m,n$ are uniformly close to the endpoints of $\sigma$. Similarly,
if $m', n'$ are the points on $(r,s,\cdot)$ closest to $p,q$ respectively, then $m',m$ are
close, and so are $n,n'$. Now consider the flowline $(p,s,\cdot)$. This flowline converges 
exponentially fast to $(p,q,\cdot)$ along the ray from $n$ to $p$. Similarly, it
converges exponentially fast to $(r,s,\cdot)$ along the ray from $m'$ to $s$. Hence away from a bounded
neighborhood of the endpoints of $\sigma$, the flowlines $(p,q,\cdot)$ and $(r,s,\cdot)$ are themselves
exponentially close. This is the kind of convexity argument which is very standard in strictly
negatively curved spaces; Mineyev's technology allows us to transplant it to $\F$.

It follows that 
$$d^\times(e^{-1}a \cdot (p,q,C_6),(p,q,C_6)) \le 2Me^{\lambda C_6}$$
By the triangle inequality, for any $n$ we estimate
$$d^\times((e^{-1}a)^n \cdot (p,q,C_6), (p,q,C_6)) \le 2nMe^{\lambda C_6}$$
and therefore
$$\tau(e^{-1}a) \le \frac {2nKMe^{\lambda C_6} + \epsilon} n$$
where $K,\epsilon$ are as in the last bullet of Theorem~\ref{Mineyev_theorem}.
The constants $K,M,\lambda,\epsilon$ all depend only on $\delta$ and $|S|$, whereas
we may choose $C_6$ as big as we like, and $n$ as big as we like. 
In particular, $\tau(e^{-1}a)$ may be made arbitrarily
small by choosing $C_6$ very big.

On the other hand, if $a$ is not equal to $e$, then $\tau(e^{-1}a) \ge C_1 > 0$
for $C_1$ depending only on $\delta$ and $|S|$. This gives a contradiction,
for sufficiently big $C_6$ (chosen depending on $C_1$).

This contradiction proves the theorem.
\end{proof}

\begin{example}\label{nonorientable_surface_example}
The upper bound $1/2$ is {\em sharp}, and is realized in
a nonabelian free group, or closed hyperbolic surface group, by \cite{Comerford_Edmunds}.
\end{example}

\begin{example}
Suppose $H = \langle h_1,h_2\rangle$ is a non-free $2$-generator subgroup of (any group) $G$. Then one can
show $\scl([h_1,h_2]) < 1/2$. However, Delzant showed (\cite{Delzant}) that in any word-hyperbolic
group there are only finitely many conjugacy classes of non-free $2$-generator subgroups.
\end{example}

If $G$ is allowed to have torsion, things become slightly more complicated.

\begin{example}\label{orbifold_example}
Let $S$ be a (surface) orbifold containing two orbifold points $p_2,p_3$ of order 
$2$ and $3$ respectively. If $\alpha$ is any embedded arc in $S$ from  $p_2$ to $p_3$, then
the boundary of a regular neighborhood of $\alpha$ has $\scl$ at most $1/12$.
For typical $S$ there are infinitely many distinct isotopy classes of such arcs $\alpha$.
\end{example}

It is straightforward to see that in any word-hyperbolic group $G$, there is a positive
first accumulation point for $\scl$ on conjugacy classes satisfying the no mirror
condition (in fact this follows directly from Theorem A and Theorem A${}^\prime$), 
but we have not been able to show that this lower bound is independent of $G$, and therefore
we pose the following

\begin{question}
Is there a universal positive constant $C$ such that the first accumulation point
for $\scl$ on conjugacy classes satisfying the no mirror condition in a hyperbolic
group is at least $C$?
\end{question}

Finally, the explicitness and universality of the constants in Theorem B, together with
Example~\ref{nonorientable_surface_example} motivates the following

\begin{question}
Can the lower bound in Theorem B be improved to $1/4$?
\end{question}

\section{Mapping class groups}\label{mapping_class_group_section}

Our theorems may be generalized to groups which are not themselves hyperbolic,
but which act suitably on $\delta$-hyperbolic
spaces. In this section and the next, we concentrate on two important examples.

We show how to adapt our Gap Theorem to the action
of the mapping class group on the complex of curves. For an introduction to this
complex and its properties, see \cite{masur-minsky} or \cite{bowditch}.
The proof follows much the same outline as the proof of Theorem A${}^\prime$. A significant
difference is that the action of $\Mod(S)$ on $\CC(S)$ is not {\em proper};
nevertheless it is {\em weakly proper} in a suitable sense 
\cite{Fujiwara_Bestvina}, and this
weak properness is enough. The technical tool we use is the {\em acylindricity}
of the action of $\Mod(S)$ on $\CC(S)$ as observed by Bowditch
\cite{bowditch}.

We denote the translation length of an element $a \in \Mod(S)$ 
on $\CC(S)$ by $\tau(a)$.
\begin{mod_grp_thm}[Mapping Class Theorem]
Let $S$ be a compact orientable surface of hyperbolic type and 
$\Mod(S)$ its mapping class group.
Then there is a positive integer $P$ depending on $S$ such that 
for any pseudo-Anosov element $a$, either there is
an $0 < n \le P$ and an element $b \in \Mod(S)$ with 
$ba^{-n} b^{-1}=a^n$, or else there exists a
homogeneous quasimorphism $h$ on $\Mod(S)$ such that 
$h(a) =1$ and $D(h) \le C(S)$,
where $C(S)$ depends only on $S$.

Moreover, let $a_i \in \Mod(S)$ be a collection of elements for which
$T = \sup_i \tau(a_i)$ is finite. Suppose that for all integers $n,m \ne 0$
and all elements $b \in \Mod(S)$ and indices $i$, there is an inequality
$$ba^nb^{-1} \ne a_i^m$$
Then there is a homogeneous quasimorphism $h$ on $\Mod(S)$
such that 
\begin{enumerate}
\item{$h(a) =1 $, and $h(a_i)=0$ for all $i$}
\item{The defect of $h$ is $\le C'(S, T)$}
\end{enumerate}
\end{mod_grp_thm}

\begin{remark}
By Thurston's classification of surface automorphisms (see e.g. \cite{Thurston_dynamics})
every element of infinite order in $\Mod(S)$ is either reducible or 
pseudo-Anosov. An element $a \in \Mod(S)$ has $\tr(a)=0$
on $\CC(S)$ if and only if it has finite order or 
it is reducible.
\end{remark}
\begin{remark}
Note the reappearance of the no mirror condition. Also note the separation theorem;
compare with the statement of Theorem A${}^\prime$.
\end{remark}
\begin{remark}
The dependence of $C(S)$ on $S$ is somewhat subtle and indirect, and does not seem to
be easy to estimate.
\end{remark}
\begin{proof}
The basic structure of the proof should be reasonably clear at this point.

Let $\CC(S)$ be the curve graph of $S$.
$\CC(S)$ is $\delta$-hyperbolic, \cite{masur-minsky}.
Any pseudo-Anosov element acts as an axial isometry.
Moreover, by Bowditch \cite{bowditch}, there exists $B(S)$ such that 
for any pseudo-Anosov element $a$, $a^B$ has a geodesic
axis in $\CC(S)$. So, in particular, 
$\tr(a) \ge 1/B$.

If there exists $n>0$ and $b \in \Mod(S)$ with 
$ba^n b^{-1}=a ^{-n}$, then one may assume
$n \le P(S)$, where $P(S)$ depends only on $S$.
This follows because the action of $\Mod(S)$ on $\CC(S)$ is
{\em acylindrical} in the sense of Bowditch, \cite{bowditch}.
Here is the precise statement of acylindricity: for any $A > 0$, there 
exists $E,M$ such that for any two points
$x,y \in \CC(S)$ with $d(x,y) \ge E$ then 
there are at most $M$ elements $b \in \Mod(S)$
such that $d(x,bx) \le A, d(y,by) \le A$.
We consider the case $A=10 \delta$.
So, in the following we assume that for all $0<n$ and $b \in \Mod(S)$,
we have 
$ba^nb^{-1} \not= a^{-n}$.

From above, one finds that there is $P'>0$ such 
that if $n \ge P'$, then (a copy of the geodesic segment) $a^{-n}$
does not appear in the $L$-neighborhood of a geodesic axis of $a^B$,
where $L$ is the constant from the section 
\ref{counting_function_section} which 
depends only on $\delta$ in this setting.
This is because, otherwise, 
 one finds that $a^m$ is conjugate to
$a^{-m}$ for some $0<m$, which is a contradiction. 

So, there exists $R(S)$ such that for any $a$ as
in the theorem, there exists $N \le R$ such that 
$a^N$ has a geodesic axis,  and 
no reverse counting happens for $a^N$.

Let $\alpha$ be a geodesic axis for $a^N$. Let $x \in \alpha$, 
and denote the subpath from $x$ to $a^N x$ as $a^N$.
Let $h=c-c_{-}$ be the quasimorphism constructed 
using the counting functions for the pairs $(a^N,1)$ and $(a^{-N},1)$.
Then, for any $n>0$, $c(a^{Nn}) \ge n$ and $c_{-}(a^{Nn})=0$ since 
there is no reverse counting.
So, $h(a^{Nn}) \ge n$ for all $n>0$.
We know 
$D(h) \le K(\delta)$ where
$K(\delta)$ depends only on $S$
by Proposition~\ref{defect_estimate}.
Therefore, given $a$, we can construct (by averaging $h$) a homogeneous quasimorphism $f$
such that $f(a)=1$ and $D(f) \le 4KR$, 
where the constant $KR$ depends only on $S$.

\medskip

The argument to prove the second part is very similar to the proof of
the second part of Theorem A${}^\prime$. Given a collection of elements $a_i$ and a
uniform upper bound $T$ on their translation lengths, if there is a translate of the
axis of $a$ which stays close to an axis of $a_i$ on a sufficiently long interval, then
by acylindricity and the pigeonhole principle, we can conclude that some power of
$a$ is conjugate to a power of $a_i$.

In slightly more detail: let $l_{a_i}$ be an axis for $a_i$, and let $v$ be an axis
for a conjugate $e = cac^{-1}$ of $a$ which stays close to $l_{a_i}$ on a sufficiently
long segment $\beta$ (we say how long in a moment). Let integers $n_j,m_j$ be chosen for which
$$|n_j\tr(a_i) - m_j\tr(a)| < C_1$$
for some fixed constant $C_1$. Given an upper bound on $T$ and $\tr(a)$, we can find
at least $C_2$ many such pairs $(n_j,m_j)$ whose absolute values are bounded by $C_3$
where $C_2$ is as big as we like, and where $C_3$ depends on $T,\tr(a),S,C_1$
(and {\em not} on the specific element $a_i$).
Let $\beta$ be longer than $C_3T,C_3\tr(a)$. Then (after possibly replacing $n_j$
by $-n_j$ for some $j$) we observe that
$$d(a_i^{n_j}e^{m_j}(p),p) \le C_1'$$
for some $C_1'$ which depends only on $C_1$ and $S$. By acylindricity and the
pigeonhole principle, there
are distinct indices (which we denote by $n_1,m_1,n_2,m_2$ respectively) for which
$$a_i^{n_1}e^{m_1} = a_i^{n_2}e^{m_2}$$
and therefore some nonzero power of $a_i$ is equal to some nonzero power of $e$, which
is itself conjugate to $a$, contrary to hypothesis.

This shows that the length $\beta$ can be bounded above in terms of $S,T,\tr(a)$.

Note that one constructs by this argument a {\em single} quasimorphism $h$
whose value grows linearly on powers of $a$, and which vanishes on all
powers of $a_i$ simultaneously. The fact that there are (typically) infinitely many $a_i$
on which $h$ vanishes is immaterial.

\medskip

We give another proof of the second claim which makes more explicit use of
acylindricity and gives the slightly better constants claimed in the statement of the theorem.

The argument is similar to the one 
to show for sufficiently large $N$, there is no 
reverse counting for $a^N$, since 
otherwise, for some $0 <n \le N$, $a^n$
is conjugate to $a^{-n}$ by the acylindricity of the action,
which is impossible. Also, see the proof of Claim 2${}^\prime$
in the proof of Theorem A${}^\prime$.

We want to show $c_{a^{\pm 2N}}(a_i)=0$ for all $i$
if $N$ is bigger than a constant depending only on 
$S$ and $T$.
Suppose $c_{a^{2N}}(a_i)>0$ for some (fixed) $i$
(the argument is precisely analogous for $c_{a^{-2N}}$).
Let $l_{a_i}$ be an axis for $a_i$, and let $l_e$ be an axis
for a conjugate $e = cac^{-1}$  of $a$ which stays close to,
namely in the $L(\delta)$-neighborhood of, $l_{a_i}$ on a 
segment $v$ such that one endpoint is $v_0$ and the other end point 
is $e^{2N}(v_0)$.

For simplicity, we assume both $l_{a_i},l_e$ are geodesics 
(replace $a,a_i$ by $a^B,a_i^B$ if necessary).
Then $v$ is $2\delta$-close to $l_{a_i}$, where $\delta$ is the 
hyperbolicity constant for $\CC(S)$ 
(here we assume $v$ is much longer than $\delta$, which follows
if we take $N$ bigger than a constant depending only on $S$).

We observe that for all $0 \le j \le N$, there exists $n_j$
such that 
$$d(a_i^{n_j} e^j(v_0),v_0), \, \, 
d(a_i^{n_j} e^{j+N}(v_0),e^N(v_0)) 
\le \tau(a_i)+4\delta \le T+4\delta$$

We set $A=T+4\delta$ and obtain corresponding constants $D(A),M(A)$ for the 
acylindricity of the action.
Now, assume $N$ is such that $d(e^N(v_0),v_0) \ge D(A)$.
In other words, $N \tau(a) \ge D(A)$ (we know $\tau(a) \ge 1/B$, so 
take $N \ge BD$).
Then, by the acylindricity, 
there are at most $M$ distinct elements in $a_i^{n_j} e^j (0 \le j \le N)$.
It follows that if $N > M$, then some non-trivial power of $a$ 
is conjugated by $c$ to a (possibly trivial) power of $a_i$, which is 
impossible.
We thereby obtain an upper bound for $N$ by $A,D(A),M(A),B$, which depend
only on $S,T$, to have $c_{a^{2N}}(a_i)>0$. 
We obtain the same upper bound from $c_{a^{-2N}}(a_i)>0$ as well. 

\end{proof}

\section{amalgamations}

In this section we adapt our theorem to the special case of an amalgamated
free product acting on its associated tree, and construct many quasimorphisms
with uniform lower bounds. One must be slightly careful:
the group $\SL(2,\Q_p)$ is an amalgam of two copies of $\SL(2,\Z_p)$. Nevertheless,
$\SL(2,\Q_p)$ is uniformly perfect, and therefore admits no nonzero homogeneous quasimorphisms
at all. As in the case of the mapping class group acting on $\CC(S)$, one must
ensure (by fiat) that the action of the amalgam on its associated tree is {\em weakly
proper}; this is guaranteed by a suitable malnormality condition.

\begin{amalgam_thm}[Amalgamation Theorem]
Let $G=A*_C B$.
Let $w$ be a word on $A,B$ which is
reduced and cyclically reduced
such that $|w| >1$, and let $\overline{w}$ denote the corresponding element in $G$.
Suppose that the double coset $C\overline{w}C$,
does not contain the element corresponding to any cyclic conjugate of $w^{-1}$.
Then there exists a homogeneous quasimorphism $h$ on $G$
such that $h(\overline{w})=1$, and the defect of $h$ is $\le 312$.
\end{amalgam_thm}
\begin{proof}
We use \cite{fujiwara:tams}.
Here is a review.
Let $\Gamma$ be the Cayley graph of $G$
for the generating set $A \cup B$.

Then a geodesic between $1$ and $a$ is exactly 
a word for $a$ reduced  as a word in $A*_C B$ (Lemma 3.1 \cite{fujiwara:tams}).
Let $w$ be a reduced word for $a$.
Assume that $w$ is cyclically reduced.
Then $w^{*}= \cdots w w w \cdots $ is an infinite reduced word, which is 
a geodesic in $\Gamma$ since $w$ is cyclically reduced.
For such $w$, let $c,c_{-1}$ be the 
counting function for $(w,1),(w^{-1},1)$.
Then, $c(a^n) \ge n$ for all $n>0$.

By the double coset condition, $c_{-1}(a^n)=0$ for 
all $n>0$, namely, no reverse counting.
 The argument is essentially in \cite{fujiwara:tams}.
Suppose $c_{-1}(a^n) > 0$. Then, there is a 
realizing path $\alpha$ from $1$ to $a^n$  which contains
$w^{-1}$ as a subword. It is shown in Lemma 3.2 \cite{fujiwara:tams} 
that one can always take a realizing path to be a geodesic.
Let $\beta=w^n$, which is a geodesic from $1$ to $a^n$.
Since both reduced words $\alpha$ and $\beta$ 
represent the same element $a^n$, by Britton's lemma,
there must be $c,c' \in C$ such that 
$c {\bar w} c'= \bar{v}$, where $v$ is some 
cyclic conjugate of $w^{-1}$.
(see section 4  \cite{fujiwara:tams} for details).
But this is prohibited.

It follows that $h(a^n) \ge n$, where 
$h=c-c_{-1}$.
It is shown in Proposition~3.1 \cite{fujiwara:tams} that 
$D(h) \le 78 $. 
When we make a quasimorphism homogeneous,
an upper bound for the defect becomes $312=78 \times 4$.
\end{proof}

\begin{remark}
There is no implied suggestion that the constant $312$ is optimal! The main
point is that it is a {\em universal} constant, which does not depend in any
way on the group $G$.
\end{remark}

\begin{remark}
If $C$ is trivial, so that $G$ is a free product, Bavard (\cite{Bavard}, Prop. 3.7.2)
obtains the (optimal) lower bound $\scl \ge 1/2$.
\end{remark}

\begin{remark}
Note that the double coset condition is the precise analogue of the no mirror condition
in the context of amalgamated free products.
\end{remark}

\begin{question}
What is the optimal constant in Theorem D?
\end{question}

The following example clarifies the geometry of the double coset condition
in the context of $3$-manifold topology.

\begin{example}
Let $M$ be a knot complement in a rational homology sphere. After
choosing an orientation for the knot and for $M$, there is a natural choice of meridian and
longitude $m,l$ on $\partial M$. The longitude is defined by the property that it
is virtually trivial in homology. So there is some surface $S$ of least
Euler characteristic in $M$ whose boundary is a multiple of the longitude.

Suppose $[\partial S] = n\cdot [l]$ in homology. Define $\chi_\Q = \frac {\chi(S)} n$.
Let $M_{p/q}$ be the result of Dehn filling $M$ along the slope $p/q$. 
Let $\gamma \subset M_{p/q}$ denote
the core of the added solid torus. As in Example~\ref{hyperbolic_surgery_example},
we have an estimate for stable commutator length
$$\scl(\gamma) \le \frac {-\chi_\Q} {2p}$$
As before, for $K$ a knot in $S^3$, this estimate becomes
$$\scl(\gamma) \le \frac {g- 1/2} {p}$$
where $g$ is the genus of a Seifert surface for $K$.

A slope $p/q$ on $\partial M$ is said to be a {\em boundary slope} if there is an
essential, oriented, proper surface $\Sigma \subset M$ such that $\partial \Sigma$
is a nonempty union of curves isotopic to the $p/q$ curve, with either orientation.

After filling $M$ to $M_{p/q}$, the manifold $M_{p/q}$ splits along
a surface $\Sigma'$ obtained by filling in the boundary components of $\Sigma$, into
two submanifolds which by abuse of notation we denote $M_A,M_B$ with fundamental
groups $A,B$. For brevity, we denote $\pi_1(\Sigma) = C$.
The core geodesic $\gamma$ intersects the two submanifolds $M_A,M_B$ efficiently,
in a collection of proper arcs which represent elements of the double coset spaces
$CAC$ and $CBC$. If the volume of $M$ (and therefore that of $M_{p/q}$) is small,
and the area of $\Sigma'$ is large, most of $M_A,M_B$ must be very thin, and have
the structure of an $I$-bundle over some subsurfaces $\Sigma_A',\Sigma_B'$ of
$\Sigma' = \partial M_A = \partial M_B$. These $I$-bundles are known as the
{\em characteristic} submanifolds of $M_A,M_B$ and we denote them by $I_A,I_B$
respectively. 

If $\Sigma'$ is connected, the boundary components of $I_A$ and $I_B$ are
contained in a single surface. Suppose that $\partial I_A$ and $\partial I_B$
are connected. The geodesic $\gamma$ is decomposed into a sequence of
geodesic segments $\gamma_j$ with $0 \le j \le i$ which are the connected components of 
$\gamma \cap (M_{p/q} - \Sigma')$. Each $\gamma_i$ is properly embedded in
$M_A$ or $M_B$, and is contained in the corresponding characteristic submanifold
$I_A$ or $I_B$. Since we are assuming $\partial I_A,\partial I_B$ are connected,
each oriented geodesic segment $\gamma_j$ can be dragged around $I_A$ (say)
into itself in such a way that the orientation is reversed at the end.
Composing these proper isotopies, we can drag all the $\gamma_j$ simultaneously
in such a way that the common endpoints of $\gamma_j$ and $\gamma_{j-1}$ agree
during the dragging, for each $j\le i$. The initial point of $\gamma_0$ and
the final point of $\gamma_i$ will not necessarily follow homotopic paths
under this dragging, and their difference is an element of $\pi_1(\Sigma')$.
Translating this into algebraic information, 
we have exhibited a conjugate of $\gamma$ as an element of the double
coset space $\pi_1(\Sigma')\gamma^{-1}\pi_1(\Sigma')$.

This example actually occurs: Nathan Dunfield \cite{Dunfield} has produced
examples of Montesinos knots with genus $1$ for which $p/1$ is an essential
slope where $p \sim 20,000$. The corresponding core geodesics $\gamma$
satisfy $\scl(\gamma) \le \frac 1 {40,000}$ and therefore by Theorem D,
$\gamma$ is conjugate into $\pi_1(\Sigma')\gamma^{-1}\pi_1(\Sigma')$.

The double coset condition shows us how to think about the geometry of
the resulting manifolds, and the way in which $\gamma$ sits inside them.
\end{example}

\end{document}